\documentstyle[twoside,amstex]{article}
\input amssym.def
\input amssym.tex
\def\bc{\begin{center}}
\def\ec{\end{center}}

\begin{document}
\thispagestyle{empty} \vspace*{3 true cm} \pagestyle{myheadings}
\markboth {\hfill {\sl Huanyin Chen}\hfill} {\hfill{\sl Separative
exchange rings in which $2$ is invertible}\hfill} \vspace*{-1.5
true cm} \bc{\large\bf SEPARATIVE EXCHANGE RINGS IN WHICH
\\
\vskip2mm $2$ IS INVERTIBLE}\ec

\vskip6mm
\bc{{\bf Huanyin Chen}\\[1mm]
Department of Mathematics, Hangzhou Normal University\\
Hangzhou 310036, People's Republic of China}\ec

\vskip4mm \begin{abstract} An exchange ring $R$ is separative
provided that for all finitely generated projective right
$R$-modules $A$ and $B$, $A\oplus A\cong A\oplus B\cong B\oplus
B\Longrightarrow A\cong B$. Let $R$ be a separative exchange ring
in which $2$ is invertible, and let $a-a^3\in R$ be regular. We
prove, in this note, that $a\in R$ is unit-regular if
$R(1-a^2)R=Rr(a)={\ell}(a)$. An element $a$ in a ring $R$ is
special clean if there exists an idempotent $e\in R$ such that
$a-e\in R$ is a unit and $aR\bigcap eR=0$. Furthermore, we prove
that $a\in R$ is special clean if $aR/ar(a^2),
R/\big(aR+r(a)\big)$ are projective, and $R(a-a^3)R=Rar(a^2)=\ell
(a^2)aR$. These also extend the corresponding results in
separative regular rings.


\end{abstract}

\vskip15mm\section{Introduction} \vskip4mm Let $R$ be a ring with
an identity. An element $a\in R$ is (unit) regular if there exists
some (unit) $x\in R$ such that $a=axa$. A ring $R$ is (unit)
regular if and only if every element in $R$ is (unit) regular. As
is well known, a ring $R$ is unit-regular if and only if every
element in $R$ is the product of an idempotent and a unit, and
that a regular ring is unit-regular if and only if it has stable
range one, i.e., $A\oplus B\cong A\oplus C\Longrightarrow B\cong
C$ for all finitely generated projective right $R$-modules $A,B$
and $C$. Following the terminology used in \cite{Ak}, we say that
$a\in R$ is special clean if there exists an idempotent $e\in R$
such that $a-e\in R$ is a unit and $aR\bigcap eR=0$. Camillo and
Khurana Theorem stated that a ring $R$ is unit-regular if and only
if every element in $R$ is special clean~\cite[Theorem 1]{CK}.

A ring $R$ is an exchange ring if for any $a\in R$ there exists an
idempotent $e\in aR$ such that $1-e\in (1-a)R$. The class of
exchange rings is very large. It includes all regular rings, all
$\pi$-regular rings, all strongly $\pi$-regular rings, all
semiperfect rings, all left or right continuous rings, all clean
rings, all unit $C^*$-algebras of real rank zero and all right
semi-artinian rings, etc. A separative ring is one whose finitely
generated projective modules satisfy the property $A\oplus A\cong
A\oplus B\cong B\oplus B\Longrightarrow A\cong B$~\cite{Ara1}. For
instances, every weakly stable exchange ring (including every
exchange ring having stable range one) and exchange ring
satisfying generalized $s$-comparability~\cite{CH}. Recently,
O'Meara proved that the condition $Rr(a)={\ell}(a)R=R(1-a)R$
characterizes elements a of a regular ring $R$ that are products
of idempotents in precisely the separative regular rings. In fact,
one of the most important open problems in regular rings is that
if every regular ring is separative~\cite{KC}.

The purpose of this note is to explore when a regular element in
separative exchange rings is unit-regular. Let $R$ be a separative
exchange ring in which $2$ is invertible, and let $a-a^3\in R$ be
regular. We prove, in Section 2, that $a\in R$ is unit-regular if
$R(1-a^2)R=Rr(a)={\ell}(a)$. In Section 3, we further prove that
$a\in R$ is special clean if $aR/ar(a^2), R/\big(aR+r(a)\big)$ are
projective, and  $R(a-a^3)R=Rar(a^2)=\ell (a^2)aR$. These also
extend ~\cite[Theorem 2.2]{CH1} and ~\cite[Theorem 15.2.1]{CH}
from regular rings to exchange ones.

Throughout, all rings are associative with an identity and all
modules are right modules. For any right modules $A$ and $B$,
$A\lesssim^{\oplus} B$ means that $A$ is isomorphic to a direct
summand of $B$. We use $A\propto B$ to stand for
$A\lesssim^{\oplus }mB$ for some $m\in {\Bbb N}$. We denote the
left (right) annihilator of $a$ in a ring $R$ by
$\ell{(a)}(r(a))$.

\vskip10mm\section{Unit-regular Elements}

\vskip4mm We start by a several lemmas which will be used in the
proofs of the main results.

\vskip4mm \hspace{-1.8em} {\bf Lemma 2.1. ~\cite[Lemma
2.1]{Ara1}}\ \ {\it Let $R$ be an exchange ring. Then the
following are equivalent:}\vspace{-.5mm}
\begin{enumerate}
\item [(1)] {\it $R$ is separative.}
\item [(2)] {\it For all finitely generated projective right
$R$-modules $A,B$ and $C$, $C\oplus A\cong C\oplus B$ with
$C\propto A,B\Longrightarrow A\cong B$.}
\end{enumerate}

\vskip4mm \hspace{-1.8em} {\bf Lemma 2.2.}\ \ {\it Let $R$ be an
exchange ring in which $2$ is invertible, and let $a-a^3\in R$ be
regular. Then $$(a-a^3)R\oplus r(a)\cong (a-a^3)R\oplus R/aR.$$}
\vskip2mm\hspace{-1.8em} {\it Proof.}\ \ Write
$a-a^3=(a-a^3)s(a-a^3)$. Then $a=axa$, where
$x=a+(1-a^2)s(1-a^2)$. Set $b=1-a$. Then $a=1-b$. Hence,
$$(1-b)b(2-b)=(1-b)b(2-b)s(1-b)b(2-b).$$ As $2\in R$ is
invertible, we get
$$(1-b)b(1-\frac{b}{2})=(1-b)b(1-\frac{b}{2})(2s)(1-b)b(1-\frac{b}{2}).$$
Hence, $b=byb$, where
$y=\frac{3}{2}-\frac{b}{2}+(1-b)(1-\frac{b}{2})(2s)(1-b)(1-\frac{b}{2})$.
Set $c=1+a$. Then
$$c(1-c)(1-\frac{c}{2})=c(1-c)(1-\frac{c}{2})(-2s)c(1-c)(1-\frac{c}{2}).$$
Thus, $c=czc$, where
$z=\frac{3}{2}-\frac{c}{2}+(1-c)(1-\frac{c}{2})(-2s)(1-c)(1-\frac{c}{2})$.
Clearly, $r(a)=(1-xa)R$ and $r(b)=(1-yb)R$. Set $e=1-xa$ and
$f=(1-yb)$. Then $eR+fR=eR+(1-e)fR$. One easily checks that
$$\begin{array}{lll}
(1-e)fa(1-e)f&=&xa(1-yb)axa(1-yb)\\
&=&xa(1-yb)a(1-yb)\\
&=&xa\big(1-y(1-a)\big)a\big(1-y(1-a)\big)\\
&=&x\big(1-y(1-a)\big)\\
&=&(1-e)f. \end{array}$$ Set $g=(1-e)fa(1-e)$. Then $g=g^2\in R$
and $eg=ge=0$. It follows that $eR+fR=(e+g)R$. Set $h=e+g$. Then
$r(a)+r(b)=eR+fR=hR$ with $h=h^2$. Further,
$$R=\big(r(a)+r(b)\big)\oplus (1-h)R=r(c)\oplus zcR.$$
Since $R$ is an exchange ring, so is $end_R\big(hR\big)$. Thus, we
can find some $C\subseteq r(c), D\subseteq zcR$ such that
$R=\big(r(a)+r(b)\big)\oplus C\oplus D$.

Construct a map $\varphi: D\to a(1-a)(1+a)D$ given by
$\varphi(d)=a(1-a)(1+a)d$ for any $d\in D$. If $a(1-a)(1+a)d=0$,
then $(1-a)(1+a)d\in r(a), a(1+a)d\in r(b)$ and $a(1-a)d\in
C\oplus D$. It is easy to verify that
$$\begin{array}{lll}
d&=&(1-a^2)d+a^2d\\
&=&(1-a)(1+a)d+\frac{1}{2}a(1+a)d-\frac{1}{2}a(1-a)d\end{array}$$
and so $d+\frac{1}{2}a(1-a)d\in (C\oplus D)\bigcap (r(a)\oplus
r(b))=0$. This implies that $d=-\frac{1}{2}a(1-a)d\in r(c)\bigcap
D\subseteq r(c)\bigcap zcR=0$. Therefore $\varphi: D\to (a-a^3)D$
is a right $R$-isomorphism. Consequently, $D\cong aD\cong
(a-a^3)D=(a-a^3)R.$ Since $R=r(a)\oplus r(b)\oplus C\oplus D,
ar(b)=r(b)$ and $aC=C$, we see that $$aR=r(b)\oplus C\oplus D.$$
Hence, $$\begin{array}{lll} R&=&r(b)\oplus C\oplus aD\oplus
(1-ax)R\\
&=&r(a)\oplus r(b)\oplus C\oplus D.
\end{array}$$ It follows that $r(a)\oplus D\cong R/aR\oplus aD$,
and therefore $(a-a^3)R\oplus r(a)\cong (a-a^3)R\oplus R/aR.$
\hfill$\Box$

\vskip4mm \hspace{-1.8em} {\bf Theorem 2.3.}\ \ {\it Let $R$ be a
separative exchange ring in which $2$ is invertible, and let
$a-a^3\in R$ be regular. If $R(1-a^2)R=Rr(a)={\ell}(a)R$, then
$a\in R$ is unit-regular.} \vskip2mm\hspace{-1.8em} {\it Proof.}\
\ Suppose $R(1-a^2)R=Rr(a)={\ell}(a)R$. Construct a map $\varphi:
(1-a^2)R\to (a-a^3)R$ given by
$\varphi\big((1-a^2)r\big)=(a-a^3)r$ for any $r\in R$. Clearly,
$\varphi$ is an $R$-epimorphism. As $a-a^3\in R$ is regular,
$(a-a^3)R$ is projective. Thus, $\varphi$ splits. So we can find a
right $R$-module $D$ such that $(1-a^2)R\cong (a-a^3)R\oplus D$.
In view of Lemma 2.2, $(a-a^3)R\oplus r(a)\cong (a-a^3)R\oplus
R/aR,$ and then $$r(a)\oplus (1-a^2)R\cong R/aR\oplus (1-a^2)R.$$

From $Rr(a)=R(1-a^2)R$, we get $(1-a^2)R\subseteq Rr(a)$. Thus we
can find some $x_1,\cdots ,x_n\in R; y_1,\cdots ,y_n\in r(a)$ such
that $1-a^2=\sum\limits_{i=1}^{n}x_iy_i$. Write
$a-a^3=(a-a^3)s(a-a^3)$ for some $s\in R$. Then one easily checks
that $$1-a^2=(1-a^2)(1+a^2sa)(1-a^2).$$ That is, $1-a^2\in R$ is
regular. So we may assume that each $x_i\in (1-a^2)R$. Construct a
map $\varphi: nr(a)\to (1-a^2)R$ given by $\varphi(z_1,\cdots
,z_n)=\sum\limits_{i=1}^{n}x_iz_i$ for any $(z_1,\cdots ,z_n)\in
nr(a)$. Then $\varphi$ is a right $R$-morphism. For any
$(1-a^2)r\in (1-a^2)R$, we have $(y_1r,\cdots ,y_nr)\in nr(a)$
such that $\varphi (y_1r,\cdots ,y_nr)=(1-a^2)r$. This shows that
$\varphi$ is a right $R$-epimorphism. As $(1-a^2)R$ is projective,
$\varphi$ splits. We conclude that $(1-a^2)R\propto r(a)$.

Write $aR=gR$ for an idempotent $g\in R$. Then $R/aR\cong (1-g)R$.
Write $(1-a^2)R=eR$ for an idempotent $e\in R$. Then $e\in
ReR={\ell}(a)R$. Thus, there are $x_1,\cdots ,x_n\in {\ell}(a)$
and $r_1,\cdots ,r_n\in R$ such that
$e=\sum\limits_{i=1}^{n}x_ir_i$. Write $g=ar$, where $r\in R$.
Then $x_ig=x_iar=0$. We infer that $x_i=x_i(1-g)$. Let
$$f=diag(1-g,\cdots ,1-g)\in M_n(R), c=e(x_1,\cdots
,x_n)f, d=f\left(
\begin{array}{c}
r_1\\
\vdots\\
r_n \end{array} \right)e.$$ Then $$e=(x_1,\cdots ,x_n)\left(
\begin{array}{c}
r_1\\
\vdots\\
r_n \end{array} \right)=cd,c=ecf,d=fde.$$ Construct an
$R$-morphism $\psi: f(nR)\to eR, fx\mapsto cx$ for any $x\in nR$.
Then $\psi$ is epimorphism. Since $eR$ is projective, we get
$eR\lesssim^{\oplus }f(nR)$. This implies that
$(1-a^2)R\lesssim^{\oplus}n\big((1-g)R\big)$, i.e.,
$(1-a^2)R\propto R/aR$. Obviously, $r(a),R/aR$ are both finitely
generated projective right $R$-modules, and then we get $r(a)\cong
R/aR$ by Lemma 2.1. Therefore $a\in R$ is unit-regular, as
asserted.\hfill $\Box$

\vskip4mm \hspace{-1.8em} {\bf Lemma 2.4.}\ \ {\it Let $R$ be a
separative exchange ring in which $2$ is invertible, and let
$a-a^3\in R$ be regular. If $(a-a^3)R\propto r(a),R/aR$, then
$a\in R$ is unit-regular.} \vskip2mm\hspace{-1.8em} {\it Proof.}\
\ Suppose that $(a-a^3)R\propto r(a),R/aR$. By virtue of Lemma
2.2, $r(a)\oplus (a-a^3)R\cong R/aR\oplus (a-a^3)R$, where $r(a),
R/aR$ and $(a-a^3)$ are all finitely generated projective right
$R$-modules. It follows by Lemma 2.1 that $r(a)\cong R/aR$, and
therefore $a\in R$ is unit-regular.\hfill $\Box$

\vskip4mm \hspace{-1.8em} {\bf Theorem 2.5.}\ \ {\it Let $R$ be a
separative exchange ring in which $2$ is invertible, and let
$a-a^3\in R$ be regular. If $R(1-a^2)R\bigcap RaR=Rr(a)\bigcap
{\ell}(a)R\bigcap RaR$, then $a\in R$ is
unit-regular.}\vskip2mm\hspace{-1.8em} {\it Proof.}\ \ Suppose
that $R(1-a^2)R\bigcap RaR=Rr(a)\bigcap {\ell}(a)R\bigcap RaR$.
Then we get $$R(a-a^3)R\subseteq R(1-a^2)R\bigcap RaR\subseteq
Rr(a),$$ and so $(a-a^3)R\subseteq Rr(a)$. As in the proof of
Theorem 2.3, $(1-a^2)aR\propto r(a)$. Likewise, we have
$(1-a^2)aR\propto R/aR$. This completes the proof by Lemma
2.4.\hfill$\Box$

\vskip4mm As an immediate consequence, we derive

\vskip4mm \hspace{-1.8em} {\bf Corollary 2.6.}\ \ {\it Let $R$ be
a separative exchange ring in which $2$ is invertible. If
$R(1-a^2)R=Rr(a)\bigcap {\ell}(a)R$, then $a\in R$ is
unit-regular.}

\vskip4mm \hspace{-1.8em} {\bf Remark 2.7.}\ \ Let $R$ be a
separative regular ring. O'Meara Theorem proved that
$R(1-a)R=Rr(a)={\ell}(a)R$ implies that $a\in R$ is a product of
idempotents~\cite[Theorem 4.1]{KC}. We ask a question: if $2$ is
invertible in $R$, whether $R(1-a^2)R=Rr(a)={\ell}(a)R$ $\big(
~\mbox{or}~ R(1-a^2)R=Rr(a)\bigcap {\ell}(a)R\big)$ imply $a\in R$
is a product of idempotents?

\vskip10mm \section{Special Clean Elements}

\vskip4mm One easily checks that every special clean element is
unit-regular. But the converse is not true. The aim of this
section is to explore conditions on a separative exchange ring
that $2$ is invertible under which a regular element is special
clean.

\vskip4mm \hspace{-1.8em} {\bf Lemma 3.1.~\cite[Lemma
15.1.1.]{CH}}\ \ {\it Let $A$ be a quasi-projective right
$R$-module. If $A_1\subseteq^{\oplus }A$ and $A=A_1+B$, then there
exists some $A_2\subseteq B$ such that $A=A_1\oplus A_2$.}

\vskip4mm \hspace{-1.8em} {\bf Lemma 3.2.}\ \ {\it Let $R$ be a
ring, and let $a\in R$ be regular. If} \vspace{-.5mm}
\begin{enumerate}
\item [(1)] {\it $aR/ar(a^2)$ is projective;} \vspace{-.5mm} \item [(2)] {\it
$ar(a^2)\cong R/\big(r(a)+aR\big)$;}\end{enumerate} {\it then
$a\in R$ is special clean.}\\ \hspace{-1.8em} {\it Proof.}\ \
Construct $X,Y,K,Z,e,h$ and $v$ as in the proof of ~\cite[Lemma
15.1.2]{CH}. Then $a-e\in U(R)$. Furthermore, we check that
$$aR\bigcap eR\subseteq (K\oplus Z)\bigcap hv(R)\subseteq (K\oplus
Z)\bigcap (X\oplus Y)=0.$$ Therefore $aR\bigcap eR=0$, which
completes the proof.\hfill$\Box$

\vskip4mm \hspace{-1.8em} {\bf Lemma 3.3.}\ \ {\it Let $R$ be an
exchange ring in which $2$ is invertible, and let $a-a^3\in R$ be
regular. If $aR/ar(a^2)$ and $R/\big(aR+r(a)\big)$ are projective,
then $$(a-a^3)R\oplus ar(a^2)\cong (a-a^3)R\oplus R/(aR+r(a)).$$}
\vskip2mm\hspace{-1.8em} {\it Proof.}\ \ Let $b=1-a$ and $c=1+a$.
As in the proof of Lemma 2.2, there are $x,y,z\in R$ such that
$a=axa, b=byb$ and $c=czc$. Furthermore, we can find some
$C\subseteq r(c), D\subseteq zcR$ such that
$R=\big(r(a)+r(b)\big)\oplus C\oplus D$. Additionally, $D\cong
aD\cong (a-a^3)D=(a-a^3)R.$

Since $aR/ar(a^2)$ is projective, the exact sequence
$$0\to ar(a^2)\hookrightarrow aR\to aR/ar(a^2)\to 0$$ splits, and then there is a right $R$-module $Z$ such that
$aR=ar(a^2)\oplus Z$. As $a\in R$ is regular, we see that $aR$ is
projective. Hence, $ar(a^2)$ is projective, and then so is
$R/\big(aR+r(a)\big)$. We infer that
$$0\to aR+r(a)\hookrightarrow R\to R/\big(aR+r(a)\big)\to 0$$
splits. Thus, there exists a right $R$-module $Y$ such that
$R=\big(aR+r(a)\big)\oplus Y$, whence $aR+r(a)$ is projective.
Write $R=aR\oplus E$. Then $$aR+r(a)=\big(aR+r(a)\big)\bigcap
\big(aR\oplus E\big)=aR\oplus \big(aR+r(a)\big)\bigcap E,$$ i.e.,
$aR\subseteq ^{\oplus }aR+r(a)$. By virtue of Lemma 3.1, we have a
right $R$-module $X\subseteq r(a)$ such that $aR+r(a)=aR\oplus X$.
This implies that $r(a)=r(a)\bigcap \big(aR\oplus X\big)=K\oplus
X$, where $K:=r(a)\bigcap aR=ar(a^2)$. Therefore
$$R=aR\oplus X\oplus Y, aR=K\oplus Z, r(a)=K\oplus
X~\mbox{and}~aR+r(a)=aR\oplus X.$$ Hence, we have
$$\begin{array}{lll}
R&=&r(a)\oplus r(b)\oplus C\oplus D\\
&=&K\oplus X\oplus r(b)\oplus C\oplus D\\
&=&aR\oplus X\oplus Y\\
&=&r(b)\oplus aC\oplus aD\oplus X\oplus Y.
\end{array}$$ Clearly, $C=aC$, and then
$D\oplus K\cong aD\oplus Y$. This implies that
$$(a-a^3)R\oplus K\cong (a-a^3)R\oplus Y,$$ as asserted.\hfill$\Box$

\vskip4mm \hspace{-1.8em} {\bf Theorem 3.4.}\ \ {\it Let $R$ be a
separative exchange ring in which $2$ is invertible, and let
$a-a^3\in R$ be regular. If}\vspace{-.5mm}
\begin{enumerate}
\item [(1)] {\it $aR/ar(a^2)$ and $R/\big(aR+r(a)\big)$ are projective;} \vspace{-.5mm} \item [(2)] {\it
$Rar(a^2)=\ell (a^2)aR=R(a-a^3)R$;}\end{enumerate} {\it then $a\in R$ is special clean.}\\
\hspace{-1.8em} {\it Proof.}\ \ Since $Rar(a^2)=R(a-a^3)R$, we get
$a(1-a^2)R\subseteq Rar(a^2)$. Thus we can find some $x_1,\cdots
,x_n\in R; y_1,\cdots ,y_n\in ar(a^2)$ such that
$a(1-a^2)=\sum\limits_{i=1}^{n}x_iy_i$. As $a(1-a^2)\in R$ is
regular, we may assume that each $x_i\in a(1-a^2)R$. Construct a
map $\varphi: n\big(ar(a^2)\big)\to a(1-a^2)R$ given by
$\varphi(z_1,\cdots ,z_n)=\sum\limits_{i=1}^{n}x_iz_i$ for any
$(z_1,\cdots ,z_n)\in n\big(ar(a^2)\big)$. Analogously to the
proof of Theorem 2.3, $\varphi$ is a right $R$-epimorphism, and
therefore $a(1-a^2)R\propto ar(a^2)$.

Write $aR+r(a)=gR$ for an idempotent $g\in R$. Then
$R/\big(aR+r(a)\big)\cong (1-g)R$. Write $(a-a^3)R=eR$. Then $e\in
ReR={\ell}(a^2)aR$. Then we have $x_1,\cdots ,x_n\in {\ell}(a^2)a$
and $r_1,\cdots ,r_n\in R$ such that
$e=\sum\limits_{i=1}^{n}x_ir_i$. Write $x_i=s_ia$ for some $s_i\in
{\ell}(a^2)$. Then $s_ia^2=0$. Write $g=ar+b$, where $r\in R,b\in
r(a)$. Then $x_ig=s_ia(ar+b)=s_ia^2r=0$. We infer that
$x_i=x_i(1-g)$. Let $$f=diag(1-g,\cdots ,1-g)\in M_n(R),
c=e(x_1,\cdots ,x_n)f, d=f\left(
\begin{array}{c}
r_1\\
\vdots\\
r_n \end{array} \right)e.$$ Similarly to the proof of Theorem 2.3,
we get an $R$-epimorphism $\psi: f(nR)\to eR, fx\mapsto cx$ for
any $x\in nR$. Therefore
$a(1-a^2)R\lesssim^{\oplus}n\big((1-g)R\big)$, i.e.,
$a(1-a^2)R\propto R/\big(aR+r(a)\big)$.

In light of Lemma 3.3,  $$(a-a^3)R\oplus ar(a^2)\cong
(a-a^3)R\oplus R/(aR+r(a)).$$ By virtue of Lemma 2.1,
$ar(a^2)\cong R/\big(r(a)+aR\big)$. This completes the proof, in
terms of Lemma 3.2.\hfill$\Box$

\vskip4mm \hspace{-1.8em} {\bf Corollary 3.5.}\ \ {\it Let $R$ be
a separative regular ring in which $2$ is invertible. Then each
$a\in R$ satisfying $Rar(a^2)=\ell (a^2)aR=R(a-a^3)R$ is special
clean.} \vskip2mm\hspace{-1.8em} {\it Proof.}\ \ Since $R$ is
regular, $aR,r(a),aR+r(a)$ and $K=ar(a^2)=aR\bigcap r(a)$ are
direct summands of $R_R$. By virtue of Lemma 3.2, we have some
$X\subseteq r(a)$ such that $aR+r(a)=aR\oplus X$. Furthermore,
there exist $Y$ and $Z$ such that $R=aR\oplus X\oplus Y$ and
$aR=K\oplus Z$. Clearly, $aR/ar(a^2)\cong Z$ and
$R/\big(aR+r(a)\big)\cong Y$. Therefore $aR/ar(a^2)$ and
$R/\big(aR+r(a)\big)$ are projective, hence the result, by Theorem
3.4.\hfill$\Box$

\vskip4mm \hspace{-1.8em} {\bf Example 3.6.}\ \ Let $V$ be an
infinite-dimensional vector space over a field ${\Bbb Z}_3$, and
let $R=M_3\big(end_{{\Bbb Z}_3}(V)\big)$. Then $R$ is weakly
stable, and so $R$ is a separative regular ring by ~\cite[Theorem
5.2.9]{CH}. Choose $a=\left(
\begin{array}{ccc}
0&1&0\\
0&0&0\\
0&0&1
\end{array}
\right)\in R$. Then $Rar(a^2)\subseteq R(a-a^3)R$. Obviously,
$a^2=a^3$, and so $R(a-a^3)R=Ra(1-a)R\subseteq Rar(a^2)$. Hence,
$Rar(a^2)=R(a-a^3)R$. Similarly, we check that $\ell
(a^2)aR=R(a-a^3)R$. Thus, $Rar(a^2)=\ell (a^2)aR=R(a-a^3)R$.
Therefore $a\in R$ is special clean, in terms of Corollary 3.5.

\vskip4mm \hspace{-1.8em} {\bf Theorem 3.7.}\ \ {\it Let $R$ be a
separative exchange ring in which $2$ is invertible, and let
$a-a^3\in R$ be regular. If}\vspace{-.5mm}
\begin{enumerate}
\item [(1)] {\it $aR/ar(a^2)$ and $R/\big(aR+r(a)\big)$ are projective;} \vspace{-.5mm} \item [(2)] {\it
$Rar(a^2)=\ell (a^2)aR=R(1-a^2)R$;}\end{enumerate} {\it then $a\in R$ is special clean.}\\
\hspace{-1.8em} {\it Proof.}\ \ By hypothesis, $Rar(a^2)=\ell
(a^2)aR=R(1-a^2)R$. Construct a map $\varphi: (1-a^2)R\to
(a-a^3)R$ given by $\varphi\big((1-a^2)r\big)=(a-a^3)r$ for any
$r\in R$. Then $\varphi$ is an $R$-epimorphism. Analogously to the
proof in Theorem 2.3, $(1-a^2)R\cong (a-a^3)R\oplus D$ for a right
$R$-module $D$. In light of Lemma 3.3, we get
$$(a-a^3)R\oplus ar(a^2)\cong (a-a^3)R\oplus R/(aR+r(a)).$$ It follows that
$$(1-a^2)R\oplus ar(a^2)\cong (1-a^2)R\oplus R/(aR+r(a)).$$
As is the proof of Theorem 3.4, we claim that $(1-a^2)R\propto
ar(a^2),R/(aR+r(a))$. Since $R$ is separative exchange ring, we
get $ar(a^2)\cong R/(aR+r(a))$. According to Lemma 3.2, $a\in R$
is special clean.\hfill$\Box$

\vskip4mm \hspace{-1.8em} {\bf Corollary 3.8.}\ \ {\it Let $R$ be
a separative regular ring in which $2$ is invertible. Then each
$a\in R$ satisfying $Rar(a^2)=\ell (a^2)aR=R(1-a^2)R$ is special
clean.} \vskip2mm\hspace{-1.8em} {\it Proof.}\ \ As in the proof
of Corollary 3.5, we prove that $aR/ar(a^2)$ and
$R/\big(aR+r(a)\big)$ are projective. Therefore the result follows
by Theorem 3.7.\hfill$\Box$

\end{document}